\documentclass{amsart}

\usepackage{graphicx}
\usepackage[english]{babel}
\usepackage{csquotes}
\usepackage{amssymb}
\usepackage{tikz-cd}

\usepackage{amsmath,amsfonts,pictex,graphicx,etex,tikz,hyperref}

\usepackage{enumerate}
\usepackage{enumitem}

\newtheorem{theorem}{Theorem}[section]
\newtheorem{lemma}[theorem]{Lemma}

\newtheorem{conjecture}[theorem]{Conjecture}

\theoremstyle{definition}
\newtheorem{Definition}[theorem]{Definition}
\newtheorem{definition}[theorem]{Definition}

\theoremstyle{remark}

\newtheorem*{Obstacle}{Obstacle}

\numberwithin{equation}{section}

\newcommand{\Z}{\mathbb{Z}}

\newcommand{\R}{\mathbb{R}}

\newcommand{\holder}{H\"{o}lder }
\newcommand{\smooth}{ C^{\infty}}

\newcommand{\doititle}{}

\def\rh{Rodriguez Hertz}

\def\arXiv#1{\@ifundefined{href}{{\mdseries\ttfamily arxiv:#1}}{\href{http://arxiv.org/pdf/#1}{{\mdseries\ttfamily arXiv:#1}}}}

\newtheorem{Structural Stability Theorem}[theorem]{Structural Stability Theorem}

\begin{document}
\date{March 2, 2016}

\title{On the work of Rodriguez Hertz on rigidity in dynamics}
\author{Ralf Spatzier}
\thanks{Article published in the Journal of Modern Dynamics, \textbf{10} (2016), 191--207}
\address{Department of Mathematics, University of Michigan, 2074 East Hall, 530 Church Street, Ann Arbor, MI 48109-1043, USA}
\email{spatzier@umich.edu}

\begin{abstract}
 This paper is a  survey about recent progress in measure rigidity and global rigidity of Anosov actions, and celebrates  the 
  profound contributions by Federico Rodriguez Hertz to rigidity in dynamical systems.     
  \end{abstract}

\maketitle

\section{Introduction}\label{section:Introduction}

Federico Rodriguez Hertz is the winner of the fifth Brin Prize in Dynamical Systems. He has made terrific contributions to smooth dynamics in general and  rigidity properties of  actions by higher rank abelian groups  in particular.  I will discuss recent progress  on rigidity  in this brief survey with particular attention on the work of Rodriguez Hertz.  Two main themes will be  measure rigidity for certain non-trivially commuting non-uniformly hyperbolic diffeomorphisms and global rigidity of Anosov diffeomorphisms with nontrivial centralizers on nilmanifolds.  

It gives me great pleasure to discuss and celebrate Rodriguez Hertz' amazing contributions.   The Brin Prize mentions four papers by F. Rodriguez Hertz in their award.  Two of these concern the topic of this survey, 
the papers \cite{KalKatRod09} on non-uniform measure rigidity and \cite{Hertz} on global rigidity of Anosov actions.  We also mention several others, in particular \cite{RH-W}  joint with Z. Wang and \cite{Brown:2015aa} joint with A. Brown and Z. Wang.  
These are the crowning achievements in the work on global rigidity of  Anosov actions on tori and nilmanifolds.  

\section{Prehistory}		 \label{section:Prehistory}

Let us set the scenery by  recalling some basic facts and definitions about Anosov (also called uniformly hyperbolic) diffeomorphisms (cf., e.g., \cite{KatokHasselBook}).

\begin{Definition}
Let $M$ be a compact manifold with some ambient Riemannian metric $\|\cdot\|$. Then, a diffeomorphism  $f: M \to M$ is called {\em Anosov} if there is a continuous splitting of the tangent bundle $TM = E^s \oplus E^u$  and constants $C>0$ and $\lambda >0$ such that for all $n >0$, $v \in E^s$  and $w \in E^u$,
\[  \|df ^n (v)  \|  \leq C e^{-n \lambda }  \|v\|  \:\:\:\text{ and }  \:\:\: \|df ^{-n} (w)  \|  \leq C e^{-n \lambda }  \|w\|   .\]
\end{Definition}

Anosov diffeomorphisms are plentiful once they exist thanks to the 

\begin{theorem}[Structural stability]
If $f: M \to M$ is Anosov, so is any $g: M \to M$ sufficiently close to f in the $C^1$-topology.  Moreover, there is a \holder homeomorphism $\phi$  homotopic and $C^0$ close  to $id _M$ such that g is \holder conjugate to $f$: $g = \phi \circ f \circ \phi ^{-1}$.
 \end{theorem}

We  call this simplest kind of rigidity  {\em structural stability} or   {\em $C^0$ local rigidity}.  The conjugacy is not differentiable in general.  Indeed, 
as is well known, any Anosov diffeomorphism $f$ has many periodic points.  Now it is easy to construct a local perturbation $g$ around a given periodic point $p$ that changes the derivative at p. 

Examples of Anosov diffeomorphisms exist on  tori $T^n = \R^n / \Z ^n$,  more generally on certain nilmanifolds $N / \Lambda$ and infra-nilmanifolds, i.e. manifolds finitely covered by a nilmanifold.  Here $N$ is a simply connected nilpotent group and $\Lambda \subset N$ a discrete cocompact group. Examples can be explicitly constructed as automorphisms of the nilpotent group $N$ which preserve $\Lambda$, and are number theoretic in nature.   Much work has been done understanding which nilpotent groups admit such lattices and automorphisms.  Classification is possible in low dimensions.    More generally, one can consider {\em affine diffeomorphisms} of infra-nilmanifolds.  These are simply Anosov diffeomorphisms which lift to affine maps on the nilpotent Lie group $N$. Here we call a map $\tilde{f}: N \to N$ {\em affine}  if $\tilde{f} = A \circ L_g$ is the composition of an automorphism $A$ of $N$ with a left translation $L_g$ of $N$. 

Smale asked  in 1967, acknowledging Anosov \cite{Smale68ICM}:

\begin{conjecture}
Is every Anosov diffeomorphism f of a compact manifold M topologically conjugate to an affine Anosov map of an infra-nilmanifold? 
\end{conjecture}

Franks and later Manning made major advances to this problem in \cite{Franks1970, Manning74}.  In particular, they showed that the conjecture holds for Anosov automorphisms of infra-nilmanifolds.   The resulting conjugacy is often called a {\em Franks-Manning conjugacy}. 

\begin{theorem}  \label{theorem:franks manning}
Let $f: N \to M$ be an Anosov diffeomorphism of an infra-nilmani\-fold M.  Then there is a \holder conjugacy $\phi: M \to M$ between f and an affine automorphism $g: M \to M$. 
Moreover,  $\phi$ can be chosen homotopic to the identity and is unique in its homotopy class.
\end{theorem}

This conjugacy  will play a special role in the global rigidity results for higher rank Anosov actions on infra-nilmanifolds.  For expanding maps, Gromov  proved the corresponding conjecture using his work on polynomial growth and earlier work of Shub \cite{GromovPolyGrowth,ShubExpandingMaps1970}.

We mention that there is a corresponding theory for Anosov flows, with similar basic structural results, albeit time changes complicate structural stability.   However, there are now many examples of Anosov flows, and there is little hope of any classification.

\section{Symmetry}			  \label{section:symmetry}

As we have just seen,  once there is one Anosov diffeomorphisms on a compact manifold $M$, there are many.  In fact, by structural stability the set of Anosov diffeomorphisms is open in the $C^1$ topology.  Thus  classification up to more than $C^0$ conjugacy is impossible.  However, as in many mathematical theories, the presence of symmetries may cause further rigidity.  For example in geometry,  if a closed manifold has    non-positive curvature or is aspherical, non-discreteness of the isometry group of the universal covering of a compact manifold forces the metric to be locally symmetric  \cite{EberleinAnnals1980,acta1982, Farb-Weinberger,Limbeek:2013aa,Limbeek:2014aa}.  

What is a symmetry of a smooth dynamical system?  The simplest answer is   a smooth map which commutes with the given dynamics.  Smale already posed the problem how generic diffeomorphisms $f: M \to M$  with ``non-trivial'' centralizers in $\text{Diff}(M )$ are.   Here we say that $f$ has trivial centralizer if $Z(f) =  \{f ^n \mid n \in \Z\}$.  Smale then conjectured that generic diffeomorphisms have ``trivial'' centralizers \cite{Smale68ICM},  say with respect to the $C^1$-topology on $\text{Diff}(M )$.  
First results were obtained by Kopell for maps on $S^1$ or $\R$ \cite{Kopell1970}.  Amazingly her basic approach is still used in modification.   Later Palis and Yoccoz showed triviality of centralizers for  an open dense set of Anosov diffeomorphisms,  for example   \cite{Palis-Yoccoz1,Palis-Yoccoz2}.  
Since, Smale's conjecture has continued to attract a lot of attention.  On the one hand, at least for the  $C^1$-topology, we now have a completely general solution \cite{BCW09}:

\begin{theorem}[Bonatti-Crovisier-Wilkinson, 2009]
Let M be a closed manifold.  Then there is a residual set of diffeomorphisms  in $\text{ Diff }^1 (M)$ with trivial centralizer. 
\end{theorem} 

On the other hand,   one cannot hope for more than  residual  results thanks  to~\cite{BCVW08}:

\begin{theorem}[Bonatti-Crovisier-Vago-Wilkinson, 2008]
Let M be any closed manifold.  Then there exist a nonempty open set  $\mathcal{O} \subset \text{ Diff }^1 (M)$ and a dense set of $\smooth$ diffeomorphisms f in $\mathcal{O}$ with uncountably many $\smooth$ diffeomorphisms in $Z(f)$.
\end{theorem}

One can hope to do better, and actually classify diffeomorphisms with nontrivial centralizers.  The last result however makes it clear that classification of such diffeomorphisms is impossible.  However,  one can pursue this  under more special conditions on the diffeomorphisms, such as uniform hyperbolicity.  We will take this up in Section \ref{section:anosov}.

Finally, one can try to understand more general types of symmetries.  For example, if a diffeomorphism is an element of a $\smooth$ action of a group with non-trivial relations, one could think of the other group elements as hidden symmetries. Case in point is the theory of smooth actions of higher rank Lie groups and their lattices which are strongly restricted. We will not pursue this here.

\section{History}\label{section:history}

There is by now a long history of rigidity results for actions of higher rank abelian groups, i.e., by $\Z ^k$ or $\R ^k$ for $k \geq 2$.  
This goes all the way back to Margulis'  superrigidity theorems  in the early 1970s \cite{Margulis-ICM1975,Margulis-book1991}  which essentially classified   finite dimensional representations of lattices $\Gamma$ in higher rank semisimple Lie groups of the non-compact type, e.g. $SL(n,\R)$ for $n \geq 3$.   Amazingly, the proof of this seemingly algebraic fact uses ideas from dynamics and ergodic theory.  The presence of higher rank abelian subgroups in $\Gamma$ consisting of semisimple elements is key.

 Later, R. J. Zimmer formulated the program of analyzing smooth actions of higher rank semisimple Lie groups and their lattices on manifolds \cite{Zimmer-ICM}. We refer to  \cite{Fisher-survey2007} for an excellent survey.   A key tool in these investigations was  Zimmer's  superrigidity theorem for measurable cocycles for finite measure  preserving actions of these groups \cite{Zimmer-Annals1980}.  It  generalizes Margulis' theorem.  For smooth volume preserving actions, it gives a measurable framing for the manifold which transforms under the group according to a linear representation (and possibly some commuting  cocycle taking values in a compact group---``compact noise'').  More generally, cocycles over finite measure preserving actions are measurably cohomologous to a linear representation--- up to compact noise.   Again, the presence of higher rank abelian groups is crucial, and has led to a separate program for understanding actions of higher rank abelian groups.  
 
Different motivation came from geometry, from the classification of Riemannian manifolds of nonpositive curvature and higher rank: every geodesic in such a space is contained in a totally geodesic flat strip.  This was accomplished in the early 1980s in a series of papers by Ballmann, Brin, Burns, Eberlein and Spatzier \cite{Ballmann,BS,MR1377265}. The flat strip condition  is the geometric equivalent to commutativity  and explicitly led to questions about commuting Anosov actions (\cite[2.4]{Burns-Katok1985}, \cite{spatzier-invitation2004}).  

Next came the investigation of local rigidity properties of algebraic actions of lattices $\Gamma$ in higher rank groups starting in the late 1980s. Typical examples are the standard linear action of $SL(n,\Z)$ on the $n$-torus $T^n$.  J. Lewis \cite{MR1058434} and S. Hurder \cite{Hurder,hurder-survey1994} proved infinitesimal and deformation rigidity, respectively, of these actions, Hurder by analyzing higher rank abelian subgroups and their dynamical properties.  These results were just shy of local rigidity.  Let us first coin the relevant notions. 

We let $\Gamma$ be a finitely generated group, and $\rho: \Gamma \to \text{ Diff}^{\infty} M $ an action of $\Gamma$.  Let $S$ be a finite set of generators. 

\begin{definition}
We call $\rho$ {\em  locally rigid}  if  any other action $\rho ^* : \Gamma \to \text{ Diff}^{\infty} M$ is $\smooth$-conjugate to $\rho$ provided that for all $s \in S$, $\rho (s)$ and $\rho ^* (s)$ are  sufficiently close to $\rho$ in the $C^1$-topology. More precisely, there will be a $\smooth$ diffeomorphism $\phi : M \to M$ such that   $ \phi \circ \rho ^* (\gamma) = \rho (\gamma) \circ \phi$  for all $\gamma \in \Gamma$. 
\end{definition}

One can consider variations of this definition, for example if  the actions are only $C^k$, or if one requires the perturbation to be close in the $C^l$ topology.  For simplicity, we will not consider these variants.  Local rigidity of algebraic actions of higher rank abelian groups  was established under weaker and weaker hyperbolicity assumptions on the action by a number of teams: A.~Katok, J.~Lewis, R.~Zimmer \cite{KL1,KL2,KL}, M. Einsiedler and T. Fisher \cite{einsiedler-fisher2007}, A. Katok and R. Spatzier \cite{KS97}, D. Damjanovic and A. Katok  \cite{DamjKatok-annals2010,DamjKatok2011}.  Most recently  K. Vinhage by himself and then with J. Wang  proved the crown jewel in this program, local rigidity of quite general algebraic partially hyperbolic higher rank abelian actions  \cite{vinhageJMD2015,Vinhage:2015aa}.

For actions of  lattices in higher rank semisimple Lie groups, local rigidity for algebraic actions  was proved in full in a series of papers by D.~Fisher, G.~Margulis, and N.~Qian, cf. \cite{FM3} for the most general results. They combined dynamical methods for higher rank abelian subgroups with cocycle superrigidity, and Kazhdan's property.   This followed earlier work by S. Hurder, A. Katok, J. Lewis,  R. Spatzier,
and  R.  Zimmer who  applied their above mentioned  local rigidity results for higher rank abelian groups to get local rigidity for higher rank lattices \cite{Hurder,hurder-survey1994,KL1,KL2,KL,KS97}.

We remark that none of  the rigidity and local rigidity results above hold when the action is a product of two Anosov diffeomorphisms $\phi _i : M_i \to M_i, i=1,2$ on $M_ 1 \times M_2$.  Indeed we can perturb the $\phi _i$ even locally.  More generally, one has to exclude the existence of a $\smooth$ rank  1 factor, and finite quotients of such. We call such actions  ``{\em irreducible}''.

Finally  we coin  similar notions, Anosov and irreducibility, for $\R^k$ actions.   In particular, such actions contain a diffeomorphism {\bf a} which acts normally hyperbolically with respect  to the orbit foliation, i.e. there is a continuous splitting of $TM$ into subbundles $E^s$, $E^u$ and $E^0$ such that $E^s$ and $E^u$ get uniformly contracted by {\bf a} in forward, respectively, backward time, and where $E^0$ is tangent to the $\R^k$ orbits.  

\begin{conjecture} [A. Katok, R. Spatzier]    \label{K-S-conjecture}
Let $\Z ^k$, $k \geq 2$, act irreducibly on a compact manifold $M$. Suppose that some $a \in \Z^k$ acts by an Anosov diffeomorphism.  Then the action is $\smooth$ conjugate to a $\Z ^k$ action by infra-nil endomorphisms on an infra-nilmanifold. 

Similarly, irreducible Anosov  $\R^k$ actions are $\smooth$ conjugate to $\R^k$ actions on $M\backslash G/\Gamma$ by left translations, where $\Gamma$ is  a lattice in  a Lie group G and $M$ is a maximal compact group which commutes with the $\R^k$ action.
\end{conjecture}   

A typical example of an $\R ^k$ Anosov action is the left action of the diagonal group on $SL(k+1, \R)/\Gamma$, where $\Gamma$ is a cocompact lattice in $SL(k+1, \R)$.  This follows easily from the basic structure theory of semisimple Lie groups.  The roots are actually the Lyapunov functionals, and the expansion is totally uniform.  

We remark that the conjecture for $\R ^k$ actions is very strong as  topological classification of Anosov flows seems impossible, and there certainly are Anosov flows on more general spaces than the $M\backslash G/\Gamma$.  Still, some progress has been made, under more stringent technical assumptions:  assume that  maximal nontrivial  intersections of stable manifolds of different elements in $\R^k$ are all one-dimensional.  These are the so-called {\em Cartan actions} treated by B. Kalinin and V. Sadovskaya in \cite{KalSad06,KalSad07} and B. Kalinin and R. Spatzier in \cite{KaSp04}.  One constructs equivariant geometric structures on these maximal intersections and uses higher rank to combine these structures  between different leaves. 
We will return to Conjecture \ref{K-S-conjecture} later in Section \ref{section:anosov}.

Actions by higher rank abelian groups also enjoy measure rigidity properties.   For algebraic higher rank hyperbolic actions, these investigations  started with H.  Furstenberg's celebrated  conjecture in the late 1960s about the classification of measures on $[0,1]$ jointly invariant under $\times 2$ and $\times 3$ \cite{Furstenberg67}, following his classification of minimal sets for the semi-group action.  After intermediate work by R. Lyons \cite{Lyons},  D. Rudolph showed that a jointly invariant probability measure $\mu$ has to be Lebesgue if $\mu$ has positive measure theoretic entropy for at least one of the two maps \cite{Rudolph}.   For homogeneous higher rank abelian hyperbolic actions this work has been vastly generalized, first by Katok and Spatzier  to higher dimensional  systems \cite{katok-spatzier-invmeasures96}, then Einsiedler, Katok, Lindenstrauss in various papers \cite{Einsiedler-Katok03,Einsiedler-Katok04,EinLInd06,EinLind10,EinLind08}.  This culminated in the recent paper by Einsiedler and Lindenstrauss about invariant measures for higher rank subgroups acting on quotients of  $S$-algebraic semisimple groups \cite{EinLInd}.  


\section{Non-uniform measure rigidity}		  \label{section:measure}

Remarkably,  hyperbolic actions by higher rank abelian groups on a compact manifold $M$  exhibit rigidity properties even when the action is not algebraic.  Naturally, the measure needs not be Haar  as the action is not algebraic.  Yet one can ask if the invariant measure is absolutely continuous on M.  This turns out to be correct, at least in low dimensions relative to the rank, and under  a genericity assumption.  To discuss the details, we first have to introduce some  technical information.  

Let   $f: M \to M$ be a $C^{1, \theta}$ diffeomorphism, i.e., a diffeomorphism whose derivatives are $\theta$-\holder.   If $\mu$ is an $f$-invariant probability measure, we have the (measurable) Oseledets decomposition of $TM$ into Lyapunov subspaces of $f$ for $\mu$.    For a $\Z ^k$ action, we have a refined partition: we get linear functionals  $\chi : \R ^k \to \R$  such that $\chi (\mathbf{a}) $ are the Lyapunov exponents of $a \in \Z^k$ w.r.t. $\mu$.  In addition, we get  a measurable decomposition $TM = \oplus E^{\chi _i}$ into {\em coarse Lyapunov subspaces}  $E^{\chi _i}$,  such that a.e. w.r.t. $\mu$, $v \in E^{\chi_i}$ has Lyapunov exponent $\chi _i (a)$ for all $ \mathbf{a} \in \Z ^k$.  

Furthermore, we call the kernel of a Lyapunov exponent a {\em Weyl chamber wall} and a connected component of the complement  of the union of the Weyl chamber walls a {\em Weyl chamber} 
 of the action.  It is easy to see that Weyl chambers are convex open subsets of $\R ^k$. 

Call $m$ hyperplanes (containing $0$) in $\R ^k$   {\em in general position} if the dimension of the intersection of any $l$ of them is the minimal possible, i.e., is equal to $\max (k-l, 0)$.  Call linear functionals {\em in general position} if their kernels are.

We can now report on the first measure rigidity results for non-algebraic actions, by B. Kalinin, A. Katok and F. Rodriguez Hertz.  This was executed over several difficult technical papers.  Here is the main result from \cite{KalKatRod09}:

\begin{theorem}[B. Kalinin, A. Katok and F. Rodriguez Hertz, 2011]
Let $\Z ^k$, $ k \geq 2$, act on a $k+1$ dimensional compact manifold by $C^{1, \theta}$ diffeomorphisms.  Let $\mu$ be an invariant measure such that the measure theoretic entropy $h_{\mu}  (\mathbf{a}) >0$ for some element $\mathbf{a} \in \Z ^k$.  Assume that the Lyapunov exponents of $\Z ^k$ w.r.t. $\mu$ are in general position.  Then $\mu$ is absolutely continuous.

A similar result holds for $\R ^k$ actions on a compact manifold of dimension $2k +1$.
\end{theorem}
 
The assumption on dimension and general position will force  the coarse Lyapunov spaces to  have dimension 1, a Cartan like condition.  This greatly simplifies the structure of the action.   This result was generalized by A.~Katok and F.~Rodriguez Hertz in \cite{Katok:2010aa}  to higher dimensional manifold assuming that coarse Lyapunov distributions are one dimensional  and the action satisfies the full entropy condition, namely elements from  each kernel of every Lyapunov exponent contribute to entropy.  Formally one assumes that the entropy functional is not differentiable at the kernel of each Lyapunov exponent.
These results are the crowning achievement of several papers, for more and more general conditions.  One earlier result for example assumed that the action was on a torus and used the semi-conjugacy to the action by automorphisms induced on homology. The result under discussion does not need 
any semi-conjugacy.  This is similar to the classification results of Cartan Anosov actions by B. Kalinin, V. Sadovskaya  and R. Spatzier  respectively which does not require information about the underlying manifolds \cite{KaSp04,KalSad06,KalSad07}.  

We remark that there are natural examples of both higher rank lattice and abelian actions on tori with one or several fixed points ``blown up'', i.e. replaced by projective spaces.  These behave much like the Anosov actions on the given torus but are on a different manifold, and are not uniformly hyperbolic \cite{KL}.

\subsection{Review of the proof for algebraic actions}  
We will recall the basic ideas from the algebraic case.    For simplicity we assume that we have an $\R ^k$ action, which we can always get by inducing from a $\Z ^k$ to an $\R ^k$ action.  This allows us to consider elements in the kernels of  Lyapunov exponents.  Here is an outline of the argument:

\begin{enumerate}
\item The  coarse Lyapunov spaces are given by eigenspaces of the automorphisms and integrate to 1-dimensional foliations  ${\mathcal W} _{\chi}$.
\item The measure $\mu$ can be decomposed into conditional measures $\mu  _{\chi}$ with respect to the ${\mathcal W} _{\chi}$.
\item  Assume there is $\mathbf{a} \in \R ^k$ which is ergodic w.r.t. $\mu$.  By ergodicity, for $\mu$ generic $x$,  suitable translates  $\mathbf{a} ^{n_l} (x)$ will approximate any point $y$ in the support of $\mu _{\chi} (x) $     on  ${\mathcal W} _{\chi}(x) $.  One can show that the push forward measures $\mathbf{a}^{n_l} _* (\mu _{\chi} (x))$ are the conditional measures $\mu _{\chi} (\mathbf{a} ^{n_i} (x))$, and limit onto $\mu _{\chi} (x) $, by a Lusin argument. It is crucial that the $\mathbf{a} ^{n_i} $ are isometries along ${\mathcal W} _{\chi}(x) $.
\item One can achieve ergodicity of the special elements $\mathbf{a}$ using the $\pi$-partition trick from ergodic theory \cite{katok-spatzier-invmeasures96}. Thus one arrives at the  
\item {\em Dichotomy}: for $\mu$-a.e. x,  either the $\mu _{\chi} (x)$ are atomic or Lebesgue.
\item If any one of the $\mu  _{\chi}  $ are Lebesgue, then $\mu$ is invariant under the linear flow along the ${\mathcal W} _{\chi _i}$ foliation, and hence is the Lebesgue measure on a subtorus.  By the dimension restriction, the torus has to be the full torus.
\item Otherwise, all the conditional measures $\mu _{\chi} $, for all $\chi$, are atomic.  One then argues that the  conditional measures on stable leaves $W^s_{\mathbf{b}}$ for $\mathbf{b}  \in \R ^k$ are atomic as well.  Now   it is well-known that  $h_{\mu}(\mathbf{b}) $ is positive precisely when the conditional measures on stable leaves are not atomic.  
\end{enumerate} 

\subsection{How to make this work for non-uniform actions}   Here are the key components that extend fairly easily from  algebraic to the non-uniformly hyperbolic actions:

\begin{enumerate}
\item  coarse Lyapunov spaces, and the corresponding measurable Oseledets decomposition;
\item stable manifolds through $\mu$-a.e. x;
\item coarse Lyapunov spaces are tangent to intersections of stable foliations for different elements in $\R ^k$;
this follows as  the  coarse Lyapunov exponents are in general position, and the dimension of $M$ is $k+1$;
\item positive measure theoretic entropy for one element b, $h _{\mu} (\mathbf{b}) >0$, implies that conditional measures are not atomic.
\end{enumerate} 

While all of these arguments generalize, with some difficulty,  one now encounters a significant   problem: 

\begin{Obstacle}  \label{disaster}
The $\mathbf{a}^{n_i}$ are not isometric along the foliation tangent to $E^{\chi}$ anymore, and may even grow subexponentially when returning to a generic leaf.  Thus convergence is problematic, and the key engine of the rigidity process is undermined.  
\end{Obstacle}

The solution required great originality and bravado:  in short,  do a measurable  time change which has good  expansion properties.
Let us give an outline of the procedure.  

\subsubsection{$\varepsilon$-Lyapunov metric}
First we make a definition, standard in Pesin theory.   
We refer to standard texts on Pesin theory for the precise definition of Pesin sets.  Essentially they are compact sets of large measure on which one has uniform estimates for the relevant derivatives.   

Fix a Riemannian metric $\langle \cdot,\cdot \rangle$ on $M$, with induced norm $\| \cdot\|$. 
We define the $\varepsilon$-{\em Lyapunov metric} as follows:  Let $u,v \in E_{\chi} (x)$, and fix $\varepsilon >0$.  Set
\[   \langle u, v \rangle _{\varepsilon} := \int _{\R ^k}  \langle (d \mathbf{s}) _* u, (d \mathbf{s}) _* w \rangle e^{-2 \chi (\mathbf{s}) - 2 \varepsilon \|\mathbf{s}\|} ds .\]

 Note that for $\mathbf{a} \in \R ^k$:
 \begin{enumerate} 
\item We have the growth estimate for the derivative $D  ^{E_{\chi}} (\mathbf{a})$ of a along $E_{\chi}$
\[ e^{\chi (\mathbf{a}) - \varepsilon \| \mathbf{a}  \|}  \leq \|D _x  ^{E_{\chi}} (\mathbf{a})  \| _{\varepsilon}  \leq   e^{\chi (\mathbf{a}) + \varepsilon \| \mathbf{a}  \|} .\]
\item The  $\varepsilon$-Lyapunov metric is $C^1$ on $\R ^k $ orbits. 
\item  On any Pesin set ${\mathcal R} ^L _{\varepsilon}$, the $\varepsilon$-Lyapunov metric is \holder continuous with exponent $\gamma > l$  and \holder constant $K= K(l,\varepsilon)$ depending only on $l$ and~$\varepsilon$. 
\item Moreover, we have the following comparisons with the given Riemannian metric  for some constants $c$ and  $c= c(l,\varepsilon)$ depending only on $l$ and $\varepsilon$:
\[  c \|u\| \leq \|u\|_{x, \varepsilon}   \leq c(l,\varepsilon) \|u\|   .\]
 \end{enumerate}
 
\subsubsection{Measurable time change}  
Find a vector $w$ normal to $\text{ker}\, \chi$ such that $\chi (w)= 1$.  
Find  a function $g: M \times  
\R^ k \to \R$ such that  
\begin{enumerate}
\item $\mid g(x,\| \mathbf{a}  \|) \mid \leq 2 \varepsilon \|\| \mathbf{a}  \|\|$, and
\item   the function $G(x, \| \mathbf{a}  \|) :=  \| \mathbf{a}  \| + g(x,\| \mathbf{a}  \|) w $  satisfies
\[  \|  D _x ^{E_{\chi}} (G(x, \| \mathbf{a}  \|)) \| _{\varepsilon}  = e^{\chi (\| \mathbf{a}  \|)} . \]
 \end{enumerate}
Then,  $G(x,\| \mathbf{a}  \|) $ is a measurable {\em time change function}   on $ \R ^k$,    
continuous, indeed \holder on Pesin sets, differentiable on $\R ^k$, and with exponent $\chi$.   

Now one can repeat the arguments from the algebraic case.

\subsection{Related works} A. Katok and F. Rodriguez Hertz established the following  related result in \cite{Katok:2013aa} which they call an {\em  arithmeticity theorem}.   

\begin{theorem} [Katok-\rh, 2016]
Suppose $k \geq 2$ and that $\Z ^k$ acts $C^r$, $r >1 + \theta$, for some $\theta >0$, on a $k+1$-dimensional closed manifold M, preserving a probability measure $\mu$.  Suppose  at least one element has positive entropy w.r.t. $\mu$, and that the Lyapunov exponents are in general position.   Then there are finitely many   disjoint measurable sets $R_i$ whose union has full measure such that:  the stabilizer subgroups of each $R_i$ are measurably conjugate to an action by toral automorphisms.  
In fact, the conjugacy is bijective a.e., differentiable on most stable manifolds, and differentiable in the sense of Whitney on a set of large measure.
\end{theorem}  

As a corollary they get that entropies are given by logarithms of units in algebraic number fields, simply because $\mu$ pushes forward to Lebesgue measure and entropy of elements in $SL(n, \Z)$ acting on tori.    Furthermore, they can also control  singularities and get strong topological corollaries.  

In \cite{KatRod10}, A. Katok and F. Rodriguez Hertz  obtained the  classification of a real analytic action $\alpha$ of $SL(n,\Z)$ on the $n$-torus $T^n$ assuming the induced map on homology is standard: on the complement of a finite set, there is a real analytic conjugacy from the  linearized action to $\alpha$ whose image has full measure.   Very recently, A. Brown, F. Rodriguez Hertz and Z. Wang  found  a  major improvement~\cite{Brown:2015aa}.  We will discuss it below in Section \ref{anosov actions lattices}.  

\section{Global rigidity of higher rank Anosov actions}\label{section:anosov}

We now turn to global rigidity properties of hyperbolic actions of higher rank abelian groups and lattices. 

\subsection{First results}    \label{anosov actions abelian early} 

The general classification of higher rank Anosov actions is an outstanding problem---cf. Conjecture \ref{K-S-conjecture}---with progress only under some fairly restrictive assumptions.  One can ask instead about classification of actions on manifolds which are known to support algebraic higher rank Anosov actions.    We already know that there are no new actions $C^1$-close to algebraic actions, thanks to the local rigidity theorems of Katok and Spatzier  \cite{KS97}.  Thus we are asking the {\em global rigidity problem}.   If at least one element acts by an Anosov diffeomorphism, 
we now have a complete classification of $\smooth$ irreducible actions of $\Z ^k$, $k \geq 2$, on infra-nilmanifolds.  We will now review the important  steps culminating in the global rigidity result by F. Rodriguez Hertz and Z. Wang, Theorem \ref{theorem:global rigidity}. 

There were early results on global rigidity of $\Z ^k$ actions on tori, e.g. by A.~Katok and J. Lewis \cite[Theorem 4.12]{KL1}, under stringent assumptions such as one-dimensionality of the coarse Lyapunov foliations,  or the classification of Cartan and conformal  Anosov actions  of $\R ^k$ and $\Z ^k$  on general manifolds, cf.  B. Kalinin and R. Spatzier \cite{KaSp04} for $k \geq3$ and B. Kalinin and  V. Sadovskaya \cite{KalSad06,KalSad07}.  Kalinin and Sadovskaya recently announced a generalization of the classification for Cartan actions to $\Z ^2$ and $\R ^2$. 
All of these results  made the crucial assumption that there are many Anosov elements.  

 F. Rodriguez Hertz proved the following   result for higher rank actions on tori in 2007  \cite{Hertz}.   It is the first global rigidity result that only assumes that there is just one Anosov element in the action.  This is the second rigidity paper mentioned explicitly in the Brin prize nomination.

\begin{theorem}  \label{RH theorem}
Suppose $A \in GL(N,\Z)$, that its characteristic polynomial is irreducible over $\Z$, and that the rank of the centralizer $Z(A)$ in $GL(N,\Z)$ is at least 2.    Suppose a  $\smooth$  action $\alpha$  of $Z(A)$ on $T^N$  has at least one Anosov element $\alpha ({\mathbf a})$, ${\mathbf a} \in Z(A)$, and that the  induced action $\rho$ on homology is the given action $\rho$  of $Z(A)$ on $\R ^N$.  Then $\alpha$ is $\smooth$ conjugate to $\rho$. 
\end{theorem}

This result follows from a more general theorem we will not discuss in detail, namely \cite[Theorem 2.1]{Hertz}. The latter  makes assumptions about properties of the linearized action which    easily follow  for $Z(A)$, $A$ as above.   

The main reason global rigidity questions on tori and nilmanifolds are more approachable than classification on general manifolds   is due  to the existence of a \holder conjugacy.  Indeed, since $\alpha ({\mathbf a})$ is Anosov, $\alpha({\mathbf a})$ is \holder conjugate to an affine Anosov diffeomorphism, necessarily $\rho ({\mathbf a})$, by the Franks-Manning conjugacy $\phi$.  It is easy to show that $\phi $ conjugates all of 
$\alpha (Z({\mathbf a}))$ to  $\rho  (Z({\mathbf a}))$.  

 In consequence,  it is good enough to   show  that the Franks-Manning  conjugacy $\phi$ is $C^1$.   Indeeed, once $\phi$ is $C^1$, one can approximate $\phi$  in the $C^1$-topology by a $\smooth$ $\phi _0$, conjugate $\alpha$ by $\phi _0$ to get a $\Z ^k$ action $C^1$ close to $\rho$.  The latter will be $\smooth$ conjugate to $\rho$ by local rigidity, finishing the proof.


The actual proof of Theorem \ref{RH theorem} is  quite complicated and introduced several new ideas and tools.  Roughly, Rodriguez Hertz analyzes  normal forms of $\alpha (Z(A))$ on stable manifolds  and their conjugacies using the linearization $\rho$.   This yields smoothness of $\phi$ on an open dense set.  He uses   holonomies along various foliations to get smoothness everywhere.  
 Here is a more detailed outline.   While many crucial arguments  are missing we hope the the reader will get a sense of the ideas.


\vspace{.3em}

{\em Uniform Growth Estimates}:  
First one can relate Lyapunov exponents  of any two \holder conjugate diffeomorphisms $f$ and $g$.  Let $\phi$ be the conjugacy and $\mu$ an $f$-invariant measure.  By Pesin theory, $\mu$ a.e., the Lyapunov spaces correspond to submanifolds which contract according to the Lyapunov exponent.  Hence the image under $\phi$ does as well, and gives a Lyapunov exponent for $g$ for $\phi _* (\mu)$.  This new Lyapunov exponent is easily estimated in terms of the \holder constant of the $\phi$.   Furthermore, if $\mu$ is a hyperbolic measure so is $\phi _* (\mu)$.

In the situation at hand, this implies that for ${\mathbf b} \in Z(A)$, if $\rho ({\mathbf b})$ is Anosov, then all Lyapunov exponents are negative.  Sometimes this yields uniform exponential contraction (cf. also \cite{schreiber98}): 
\vspace{.3em}

{\em Suppose $f: M \to M$ is a $C^1$ diffeomorphism and $E$ a continuous invariant bundle for f.  If all Lyapunov exponents for E for all invariant measures for f on M are negative, then f uniformly contracts E.}
\vspace{.3em}

It is important here that $E$ is continuous and  all Lyapunov exponent for $E$ are negative,  Indeed, there are $C^r$ diffeomorphisms on tori \holder conjugate to linear ones which are not Anosov, for some $1< r < 3$ \cite{Gogolev2010}.

In light of uniform exponential contraction by Lyapunov exponents, 
 it is important to see when Lyapunov exponents are negative.   For hyperbolic measures $\mu$, Lyapunov exponents of $\mu$ can be approximated by Lyapunov exponents at periodic points by work of W. Sun and Z. Wang \cite{Sun-Wang2010} following \cite{Katok1980}.   

\rh~applies this to show that a continuous splitting of the stable into slow and fast stable subspaces is $C^s$ along the stable foliation. 
\vspace{.3em}

{\em Linearization and Conjugacies}:  Here \rh~plays a subtle variation  of a well-known scheme.  First he shows that one can linearize a $\Z^k$ action at a fixed point with a contraction.  This is a normal forms argument, and yields a local $C^{1 + \varepsilon}$ conjugacy.  Then he analyzes conjugacies between linear actions at fixed points assuming one of the linear actions is {\em rich}, i.e. that one can separate each coarse Lyapunov space for one linear action as a stable space for some element.  Richness is easy to  show for the linearization $\rho (Z(A))$.  He concludes that local conjugacies between a linear action and a rich linear action respect the splittings into coarse Lyapunov spaces.

In addition, the assumptions on the linearization $\rho$ imply that the eigenvalues $\lambda ({\mathbf b}), {\mathbf b} \in Z(A)$ are dense in $\R$ for any Lyapunov exponent $\lambda$.  If these coarse Lyapunov spaces  are 1-dimensional real of complex eigenspaces, this easily implies that any local conjugacy has to be smooth away from the periodic point.

The idea now is to compare local linearization on stable manifolds at periodic points with the global linearization $\rho$.





\vspace{.3em}

 {\em  The finale}: Consider a periodic point $p$.  From the above,we see that $\phi$ is $C^1$ away from $p$ on coarse Lyapunov subspaces at $p$ away from $p$.   To prove $C^1$ at $p$, lift the action and $\phi$ to $\R ^n$. Consider  $p +  z, 0 \neq z \in \Z ^n$,  different   lifts.  Then $\phi$ is $C^1$ on a coarse Lyapunov space of $p + z$ away from $p + z$.   Now consider holonomies along  transversal unstable manifolds from the coarse Lyapunov spaces in question through $p+z$ to that of $p$, both for $\alpha$ and for $\rho$.  The conjugacy naturally intertwines the holonomies.  From dynamics, the holonomy map for $\alpha$ is absolutely continuous with non-zero Jacobian. For $\rho$ is  of course just linear.  Moreover, on the linear side, the image of p under  the holonomy is  never $p + z$ since the unstable manifold for $\rho(b)$ never has any rational slopes.    
 
 From dynamics again, the Lyapunov exponents of $\alpha$ and $\rho$ coincide  since $\phi$ now has non-zero Jacobian along the coarse Lyapunov subspace.  So suppose $\rho ({\mathbf b})$ is Anosov. By the above the Lyapunov exponents at periodic points for $\rho({\mathbf b}) $ and $\alpha ({\mathbf b})$ coincide.  By approximation of general Lyapunov exponents by periodic ones, no Lyapunov exponents are 0, and hence $\alpha ({\mathbf b})$ is also Anosov.  Then from dynamics, under the special assumptions, stable foliations are $C^1$, hence holonomies and then $\phi$ are $C^1$ near $p$.


 
 

\subsection{Global rigidity  on tori and nilmanifolds }    \label{anosov actions abelian } 
The crowning achievement  here is the following result which completely and totally resolves the global rigidity problem on infra-nilmanifolds \cite{RH-W}:

\begin{theorem} [F. \rh-Z. Wang]			 \label{theorem:global rigidity}
Suppose $k \geq 2$ and that $\alpha$ is a   $\smooth$ $\Z ^k$  action on an infra-nilmanifold  manifold M with at least one Anosov element.   Suppose the linearization $\rho$ of $\alpha$ does not have rank 1 factors.  Then $\alpha$ is conjugate to $\rho$ by a $\smooth$ diffeomorphism.  
\end{theorem}  

This follows prior work by D. Fisher, B. Kalinin and R. Spatzier \cite{FKS} where this result was proved under the additional assumption   that $\alpha$ has many Anosov diffeomorphisms.  More precisely, they assume that there is at least one Anosov element in each  \hypertarget{Weylch}{Weyl chamber}.  As we will see, both this result and the methods used will play an important role in the proof of Theorem \ref{theorem:global rigidity}.   For one, the basic strategy of \rh~and Wang is to show that every Weyl chamber has an Anosov element, thus reducing their result to the earlier one in \cite{FKS}.   However, much more finesse is needed. 


Let us first understand the basic issue.  It  is not difficult to see that all elements in an open Weyl chamber $\mathcal{C}$ are Anosov provided one of them is.  On a Weyl chamber wall $\text{ker } \chi$ however, $\chi$ is  0, and thus a wall cannot contain any Anosov elements.  So one needs to prove that there are again Anosov elements $\alpha ({\mathbf b})$ after crossing $\text{ker } \chi$. 
 Note that  the stable space for $\alpha ({\mathbf b})$ will change for sure.  Indeed, suppose that you are at a point with  Oseledets  decomposition for $\alpha$: for  ${\mathbf a} \in \mathcal{C}$,  the stable space decomposes as $E^s _{\alpha ({\mathbf a})} = E^{ss} _{\alpha ({\mathbf a})} + E^{\chi} $ where $E^{ss} _{\alpha ({\mathbf a})}$ is the sum of Lyapunov spaces $E^{\lambda}$ for suitable $\lambda$ not positively proportional to $\chi$.   Then $E^{ss} _{\alpha({\mathbf a})}$ remains stable and $E^{\chi}$ becomes unstable when we cross  the Weyl chamber wall  $\text{ker } \chi$.   If $-\chi$ is also a Lyapunov exponent, then $E^{-\chi}$ will become stable crossing the wall.  Thus it is easy to understand the changes to the stable manifold  if Lyapunov theory applies.  The difficulty in general 
is that the Oseledets decomposition is only measurable and only gives non-uniform estimates. 

In the end, \rh~and Wang use the following criterion of R. Ma\~{n}e \cite{Mane77}:

{\em
A diffeomorphism  $f: M \to M$ of a compact manifold M is Anosov  if and only if the dimensions of the stable manifolds at all periodic points are the same and for all non-zero tangent vectors v:
$$\sup _{i \in \Z}\|d f^i (v) \| = \infty .$$
}

Thanks to the Franks-Manning conjugacy $\phi$, it is easy to verify that all stable manifolds at periodic points have the same dimension.  The second condition is easily checked for tangent vectors at points with  Oseledets decomposition into Lyapunov subspaces.  The difficulty is to show this is true  for all non-zero vectors.  

To resolve this,  \rh~ and Wang study regularity properties of the Franks-Manning conjugacy $\phi$.  More precisely, consider ${\mathbf a} \in \mathcal{C}$ close to the Weyl chamber wall $\text{ker } \chi$, and decompose  the target group $G$ into a fast stable $g ^{ss} _{{\mathbf a}}$,  slow  stable $g ^0 _{{\mathbf a}}$ and unstable subspaces $g ^u _{{\mathbf a}}$ where  $g^0 _{{\mathbf a}}$ consists of eigenspaces of weights proportional to $\chi$.   
Then they project the conjugacy (or more precisely its lift to the universal cover)  to $g^0 _{{\mathbf a}}$, or more precisely the corresponding subgroup.

\begin{lemma}  \label{regularity}
The projection $\phi ^0 _a$ is smooth 
 when restricted to any stable manifold $W^s _{\alpha ({\mathbf a})}$, with \holder transverse dependence.
 \end{lemma} 
 
 The point here is that transformation properties for elements close to $\text{ker } \chi$ coming from the algebraic side only produce very slow growth.  This has to be compared with the expansion/contraction on the non-algebraic side.  In summary, we can have small exponential expansion.  However, as in \cite{FKS}, \rh~and Wang use  exponential mixing for algebraic $\Z^ k$ actions on nilmanifolds.  The idea is that the exponential mixing of the group action overcomes slow exponential growth, and allows to conclude that certain sums converge, at least as distributions dual to \holder functions.  Then ideas from harmonic analysis, especially analysis of  wave front sets,  finish the proof.  Indeed, distributions that have partials of all orders dual to \holder functions along  transversal sets of \holder foliations are $\smooth$ functions.
 
 More precisely, $\phi ^0 _{{\mathbf a}}$  can be written as a series which will always converge to a \holder function.  Differentiating term by term,    the contributions to the derivatives coming from  contracting along and moving back along the stable manifold $W^s _{\alpha ({\mathbf a})}$ work against each other.   In rank 1 in particular, i.e. for a single Anosov diffeomorphism, we cannot conclude anything as expansion and contraction roughly cancel each other out, and only yield uncontrollable growth properties.  In higher rank however,  derivatives of the summands under elements close to  $\text{ker } \chi$ will only expand with a slow exponential rate.  Then one can make the series of derivatives converge as   distributions  using exponential decay of matrix coefficients.

 Exponential mixing   was established  for algebraic $\Z^ k$ actions on nilmanifolds by A. Gorodnik and R. Spatzier in \cite{GoroSpa2}.   For actions on tori, exponential mixing   is much easier,  and follows from more standard Fourier series argument and  Katznelson's lemma.   We refer to D. Fisher, B. Kalinin and R. Spatzier \cite{FKS} who followed D. Lind's approach for  a single automorphism in \cite{Lind82}.


Next, \rh~and Wang show that $E^{ss} _{\alpha ({\mathbf a})}$ is more than measurable:  For simplicity assume that $M$ is a torus where projections to subgroups are just linear  projections: let $V$ be the coarse Lyapunov, determined by $\chi$  and as above project the conjugacy $\phi$ to $V$ to get  a map $\phi ^0 _{{\mathbf a}}$.  Then $\phi ^0 _{{\mathbf a}}$ is $\smooth$ along all stable manifolds with \holder transverse dependence by Lemma \ref{regularity}.

The next lemma shows that $\phi ^0 _{{\mathbf a}}$ is also everywhere regular.  Since the kernel of $\phi ^0 _{\mathbf a}$ along $ W^s _{\alpha ({\mathbf a})}$ is exactly  $E^{ss} _{\alpha ({\mathbf a})}$, we get the following lemma.

\begin{lemma}  \
 $E^{ss} _{\alpha ({\mathbf a})}$ is  tangent to a smooth subfoliation of $ W^s _{\alpha ({\mathbf a})}$.  Moreover,  $\phi ^0 _{{\mathbf a}}$ has constant rank along all stable manifolds.
\end{lemma} 

This lemma allows to  verify the second condition  of Ma\~{n}e's criterion for being Anosov.  Indeed, let {\bf b} be an element just across the Weyl chamber wall  $\text{ker } \chi$.  W.l.o.g. we may assume that $v \neq 0$ is a stable vector.   If $v \in E^{ss} _{\alpha ({\mathbf a})}$, $\alpha ({\mathbf  b})$  expands $v$ in negative time.  Otherwise, the derivative  along the stable manifolds $D ((\phi ^0 _{{\mathbf a}})_{W^s}) (v) \neq 0$, and thus grows arbitrarily large under $\rho ({\mathbf b})$.  Since 
\[ 	D \rho({\mathbf b} ^n)	\left( D ((\phi ^0 _{\alpha ({\mathbf a})})_{W^s}) (v)  \right)		= 	 D ((\phi ^0 _{\alpha ({\mathbf a})})_{W^s}) \left(D \alpha ({\mathbf b} ^n) (v) \right)	\]
and the stable derivative $D ((\phi ^0 _{\alpha ({\mathbf a})})_{W^s})$ is uniformly bounded  over $M$ by the last lemma,  $\left(D \alpha ({\mathbf b} ^n) (v) \right)$ has to grow arbitrarily large as well.  Thus the hypotheses for Ma\~{n}e's criterion are satisfied.


The proof of the lemma itself is difficult, highly original and frankly stunning.  First note that the second claim implies the first as $E^{ss} _{\alpha ( \mathbf a)}$  is the kernel of $D \phi ^0 _{{\mathbf a}}$.  So assume $(\phi ^0 _{{\mathbf a}}) _{W^s}$ does not have constant rank.  Then the singular set is closed and  $\Z ^k$-invariant, and hence supports an ergodic $\Z ^k$-invariant measure $\nu$ (which could be a periodic orbit).   The idea now is to analyze the Jacobian in two ways.  
On the singular set, the ratio of stable volume of small neighborhoods and the volume of their image under the conjugacy tends to 0.  On the other hand, the pull back measure $\mu = (\phi ^{-1}) _* (\lambda $) of Haar measure $\lambda$ is absolutely continuous.  Hence where defined, the Jacobian of $\phi ^s$ is not 0. In fact, there is a positive lower bound which suffices  to  yield a contradiction.

\subsection{Applications to actions by lattices}  \label{anosov actions lattices} 

A. Brown, F. \rh~and Z.~Wang recently applied Theorem \ref{anosov actions abelian }  to get the following \cite{Brown:2015aa}:

 \begin{theorem}[Brown-\rh -Wang, 2015]   \label{anosov lattices nilmanifolds}
 Suppose $\Gamma$ is a lattice in a semisimple real Lie group all of whose factors have  real rank at least  2.  Suppose $\Gamma$ acts on a  nilmanifold by $\alpha$,  that some element acts by an  Anosov diffeomorphism, and that  $\alpha$ can be lifted to the universal cover.  Then $\alpha$ is $\smooth$ conjugate to the linearization on a finite index subgroup.
 \end{theorem}

This is a significant breakthrough in the Zimmer program of studying actions of higher rank groups and lattices.  Strikingly, Theorem \ref{anosov lattices nilmanifolds} makes no assumptions on  invariant measures, unlike in most results in the Zimmer program.

\end{document}